# CONSISTENT NEWTON-RAPHSON VS. FIXED-POINT FOR VARIATIONAL MULTISCALE FORMULATIONS FOR INCOMPRESSIBLE NAVIER–STOKES

D. Z. TURNER, K. B. NAKSHATRALA, AND K. D. HJELMSTAD

ABSTRACT. The following paper compares a consistent Newton-Raphson and fixed-point iteration based solution strategy for a variational multiscale finite element formulation for incompressible Navier–Stokes. The main contributions of this work include a consistent linearization of the Navier–Stokes equations, which provides an avenue for advanced algorithms that require origins in a consistent method. We also present a comparison between formulations that differ only in their linearization, but maintain all other equivalences. Using the variational multiscale concept, we construct a stabilized formulation (that may be considered an extension of the MINI element to nonlinear Navier-Stokes). We then linearize the problem using fixed-point iteration and by deriving a consistent tangent matrix for the update equation to obtain the solution via Newton-Raphson iterations. We show that the consistent formulation converges in fewer iterations, as expected, for several test problems. We also show that the consistent formulation converges for problems for which fixed-point iteration diverges. We present the results of both methods for problems of Reynold's number up to 5000.

## 1. INTRODUCTION

Within the context of variational multiscale formulations for incompressible Navier-Stokes, most methods employ a fixed-point iteration solution strategy. Using fixed-point iteration provides a simple way to linearize the nonlinear convective term of the Navier-Stokes equations and also is relatively easy to implement. An alternative approach employs Newton-Raphson iterations to obtain the solution. As opposed to fixed-point iteration, for Newton-Raphson, one must construct the tangent matrix by taking derivatives of the weak form, including the nonlinear convective term. Although there are a number of drawbacks to using Newton-Raphson, the primary advantages include quadratic convergence and a framework on which to implement advanced algorithms such as continuation methods [1]. The development of homotopic or arc-length type methods shows potential for capturing high Reynolds flows otherwise unattainable with existing methods [2].

Obtaining numerical solutions to the Navier-Stokes equations intrinsically involves a number of challenges. The most formidable include instabilities in the pressure that arise due to the underlying Stokes problem, potential instabilities in the velocity for high Reynolds flows, and dealing the with nonlinearity of the convective term. In this work, we use a mixed finite element formulation that treats the velocity and pressure as independent variables. Although mixed formulations require more degrees of freedom and result in discrete systems that are indefinite (they have positive and negative eigenvalues), they pose the benefits of practicality and robustness in extreme cases [3].

The numerical instabilities that arise stem from the saddle point nature of the underlying Stokes problem. This feature is mathematically explained by the Ladyzhenskaya-Babuška-Brezzi (LBB) stability condition [4]. It is well known that the standard Galerkin formulation does not yield stable results, but a number of stabilized formulations have been developed to combat this deficiency. For a comprehensive overview of these methods see [5, 6, 7, 2] and the references therein. We generate both the fixed-point iteration and Newton–Raphson strategies in the context of a variational multiscale framework. The variational multiscale concept [8] acknowledges a duality of scales inherent in the physics of the problem: a coarse or resolvable scale, and a fine-scale behavior that cannot be accurately captured by the mesh. To remedy instability, we model the fine-scale behavior by splitting the problem into two subproblems a coarse-scale subproblem and a fine-scale subproblem. We then construct our formulations in a way that models the fine-scale behavior using bubble functions. For the fixed-point iteration strategy, we solve the fine-scale subproblem in terms of the coarse-scale variables and substitute the result into the coarse-scale subproblem. For the consistent Newton–Raphson technique, we treat the two subproblems as separate nonlinear residuals and solve the system using an iterative update approach.

The primary contributions of this work include the development of a consistent, stabilized, Newton-Raphson type formulation for incompressible Navier-Stokes and a direct comparison with a fixed-point iteration strategy. This work also lays the groundwork for advanced solution algorithms such as arc-length or homotopy methods that require origins in a consistent linearization. The work presented herein is the extension of work presented in [9, 10].

In the following sections, we introduce the variational multiscale concept, followed by a formulation of the fixed-point iteration strategy. We then propose a Newton-Raphson solution technique and derive the associated tangent stiffness matrix. The results presented focus on comparing the Newton-Raphson technique to fixed-point iteration for a number of example problems. For all of the example problems, we employ linear triangular elements. In [11], we show that the consistent formulation is unstable for linear quadrilateral and hexahedral elements. We conclude with a summary of the main findings of this work.

## 2. GOVERNING EQUATIONS

Let $\Omega$ be a bounded open domain, and $\Gamma$ be its boundary, which is assumed to be piecewise smooth. Mathematically, $\Gamma$ is defined as $\Gamma := \bar{\Omega} - \Omega$, where $\bar{\Omega}$ is the closure of $\Omega$. Let the velocity vector field be denoted by $\boldsymbol{v} : \Omega \to \mathbb{R}^{nd}$, where "$nd$" is the number of spatial dimensions. Let the (kinematic) pressure field be denoted by $p : \Omega \to \mathbb{R}$. As usual, $\Gamma$ is divided into two parts, denoted by $\Gamma^{\boldsymbol{v}}$ and $\Gamma^{\boldsymbol{t}}$, such that $\Gamma^{\boldsymbol{v}} \cap \Gamma^{\boldsymbol{t}} = \emptyset$ and $\Gamma^{\boldsymbol{v}} \cup \Gamma^{\boldsymbol{t}} = \Gamma$. $\Gamma^{\boldsymbol{v}}$ is the part of the boundary on which velocity is prescribed, and $\Gamma^{\boldsymbol{t}}$ is part of the boundary on which traction is prescribed. The governing equations for incompressible Navier–Stokes flow can be written as

$$\boldsymbol{v} \cdot \nabla \boldsymbol{v} - 2\nu \nabla^2 \boldsymbol{v} + \nabla p = \boldsymbol{b} \quad \text{in} \quad \Omega \tag{1}$$

$$\nabla \cdot \boldsymbol{v} = 0 \quad \text{in} \quad \Omega \tag{2}$$

$$\boldsymbol{v} = \boldsymbol{v}^{\mathrm{p}} \quad \text{on} \quad \Gamma^{\boldsymbol{v}} \tag{3}$$

$$-p\boldsymbol{n} + \nu(\boldsymbol{n} \cdot \nabla)\boldsymbol{v} = \boldsymbol{t^n} \quad \text{on} \quad \Gamma^{\boldsymbol{t}} \tag{4}$$

where $\boldsymbol{v}$ is the velocity, $p$ is the kinematic pressure (pressure divided by density), $\nabla$ is the gradient operator, $\nabla^2$ is the Laplacian operator, $\boldsymbol{b}$ is the body force, $\nu > 0$ is the kinematic viscosity, $\boldsymbol{v}^{\mathrm{p}}$ is the prescribed velocity vector field, $\boldsymbol{t^n}$ is the prescribed traction, and $\boldsymbol{n}$ is the unit outward normal vector to $\Gamma$. Equation (1) represents the balance of linear momentum, and equation (2) represents the continuity equation for an incompressible continuum. Equations (3) and (4) are the Dirichlet and Neumann boundary conditions, respectively.

2.1. **Classical mixed formulation.** First, let us define the relevant function spaces for the velocity, $\boldsymbol{v}(\boldsymbol{x})$ and the weighting function associated with velocity, which is denoted by $\boldsymbol{w}(\boldsymbol{x})$. They are, respectively, defined as

$$\mathcal{V} := \{\boldsymbol{v} \mid \boldsymbol{v} \in (H^1(\Omega))^{nd}, \boldsymbol{v} = \boldsymbol{v}^{\mathrm{p}} \text{ on } \Gamma^v\} \tag{5}$$

$$\mathcal{W} := \{\boldsymbol{w} \mid \boldsymbol{w} \in (H^1(\Omega))^{nd}, \boldsymbol{w} = \boldsymbol{0} \text{ on } \Gamma^v\} \tag{6}$$

where $H^1(\Omega)$ is a standard Sobolev space [4]. In the classical mixed formulation, the function space for the pressure $p(\boldsymbol{x})$ and its corresponding weighting function $q(\boldsymbol{x})$ are given by

$$\mathcal{P} := \{p \mid p \in L^2(\Omega)\} \tag{7}$$

where $L^2(\Omega)$ is the space of square-integrable functions on the domain $\Omega$. In the fixed-point iteration based variational multiscale method, the function space for $p(\boldsymbol{x})$ and $q(\boldsymbol{x})$ will be

$$\underline{\mathcal{P}} := \{p \mid p \in H^1(\Omega)\} \tag{8}$$

For further details on function spaces refer to Brezzi and Fortin [4].

**Remark 1.** *When Dirichlet boundary conditions are imposed everywhere on the boundary, that is $\Gamma^t = \emptyset$, the pressure can be determined only up to an arbitrary constant. In order to define the pressure field uniquely, it is common to prescribe the average value of pressure,*

$$\int_\Omega p \, \mathrm{d}\Omega = p_0 \tag{9}$$

*where $p_0$ is arbitrarily chosen (and can be zero). Then, the appropriate function spaces for the pressure that should be used instead of $\mathcal{P}$ (defined in equation (7)) is*

$$\mathcal{P}_0 := \{p \mid p \in L^2(\Omega), \int_\Omega p \, \mathrm{d}\Omega = 0\} \tag{10}$$

*Another way to define the pressure uniquely is to prescribe the value of the pressure at a point, which is computationally the most convenient.*

For convenience, we define the $L^2$ inner-product over a spatial domain, $K$, as

$$(\boldsymbol{a}, \boldsymbol{b})_K = \int_K \boldsymbol{a} \cdot \boldsymbol{b} \, \mathrm{d}K \tag{11}$$

The subscript $K$ will be dropped if $K$ is the whole of $\Omega$, that is $K = \Omega$. The classical mixed formulation (which is based on the Galerkin principle) for the incompressible Navier-Stokes equations can be written as: Find $\boldsymbol{v}(\boldsymbol{x}) \in \mathcal{V}$ and $p(\boldsymbol{x}) \in \mathcal{P}$ such that

$$(\boldsymbol{w}, \boldsymbol{v} \cdot \nabla \boldsymbol{v}) + (\nabla \boldsymbol{w}, 2\nu \nabla \boldsymbol{v}) - (\nabla \cdot \boldsymbol{w}, p) = (\boldsymbol{w}, \boldsymbol{b}) + (\boldsymbol{w}, \boldsymbol{h})_{\Gamma_t} \quad \forall \, \boldsymbol{w} \in \mathcal{W} \tag{12}$$

$$(q, \nabla \cdot \boldsymbol{v}) = 0 \quad \forall \, q \in \mathcal{P} \tag{13}$$

To determine a numerical solution using the finite element method, one first chooses the approximating finite element spaces, which (for a conforming formulation) will be finite dimensional subspaces of the underlying function spaces of the weak formulation. Let the finite element function spaces for the velocity, the weighting function associated with the velocity, and the pressure be denoted by $\mathcal{V}^h \subseteq \mathcal{V}$, $\mathcal{W}^h \subseteq \mathcal{W}$, and $\mathcal{P}^h \subseteq \mathcal{P}$ respectively. The finite element formulation of the classical mixed formulation then reads: Find $\boldsymbol{v}^h(\boldsymbol{x}) \in \mathcal{V}^h$ and $p^h(\boldsymbol{x}) \in \mathcal{P}^h$ such that

$$(\boldsymbol{w}^h, \boldsymbol{v}^h \cdot \nabla \boldsymbol{v}^h) + (\nabla \boldsymbol{w}^h, 2\nu \nabla \boldsymbol{v}^h) - (\nabla \cdot \boldsymbol{w}^h, p^h) = (\boldsymbol{w}^h, \boldsymbol{b}) + (\boldsymbol{w}^h, \boldsymbol{h})_{\Gamma_t} \, \forall \, \boldsymbol{w}^h \in \mathcal{W}^h \tag{14}$$

$$(q^h, \nabla \cdot \boldsymbol{v}^h) = 0 \, \forall \, q^h \in \mathcal{P}^h \tag{15}$$

For mixed formulations, the inclusions $\mathcal{V}^h \subseteq \mathcal{V}$, $\mathcal{W}^h \subseteq \mathcal{W}$, and $\mathcal{P}^h \subseteq \mathcal{P}$ are themselves not sufficient to produce stable results, and additional conditions must be met by these finite element spaces to obtain meaningful numerical results. A systematic study of these types of conditions on function spaces to obtain stable numerical results is the main theme of *mixed finite elements*. One of the main conditions to be met is the LBB *inf-sup* stability condition. For further details, see [4, 2].

## 3. VARIATIONAL MULTISCALE FRAMEWORK

To remedy the inadequacies of the standard Galerkin formulation as presented above, both the fixed-point iteration based and the Newton–Raphson based solution strategies incorporate a variational multiscale framework for stabilization. The variational multiscale concept, which stems from the pioneering work by Hughes [8], decomposes the underlying fields into coarse or resolvable scales and subgrid or unresolvable scales. This decomposition provides a systematic way to develop stabilized methods.

3.1. **Multiscale decomposition.** Let us divide the domain $\Omega$ into $N$ non-overlapping subdomains $\Omega^e$ (which in the finite element context will be elements) such that

$$\Omega = \bigcup_{e=1}^{N} \Omega^e \tag{16}$$

The boundary of the element $\Omega^e$ is denoted by $\Gamma^e$. We decompose the velocity field $\boldsymbol{v}(\boldsymbol{x})$ into coarse-scale and fine-scale components, indicated as $\bar{\boldsymbol{v}}(\boldsymbol{x})$ and $\boldsymbol{v}'(\boldsymbol{x})$, respectively.

$$\boldsymbol{v}(\boldsymbol{x}) = \bar{\boldsymbol{v}}(\boldsymbol{x}) + \boldsymbol{v}'(\boldsymbol{x}) \tag{17}$$

Likewise, we decompose the weighting function $\boldsymbol{w}(\boldsymbol{x})$ into coarse-scale $\bar{\boldsymbol{w}}(\boldsymbol{x})$ and fine-scale $\boldsymbol{w}'(\boldsymbol{x})$ components.

$$\boldsymbol{w}(\boldsymbol{x}) = \bar{\boldsymbol{w}}(\boldsymbol{x}) + \boldsymbol{w}'(\boldsymbol{x}) \tag{18}$$

We further make an assumption that the fine-scale components vanish along each element boundary.

$$\boldsymbol{v}'(\boldsymbol{x}) = \boldsymbol{w}'(\boldsymbol{x}) = \boldsymbol{0} \quad \text{on} \quad \Gamma^e \, ; \, e = 1, \ldots, N \tag{19}$$

Let $\bar{\mathcal{V}}$ be the function space for the coarse-scale component of the velocity $\bar{\boldsymbol{v}}$, and $\bar{\mathcal{W}}$ be the function space for $\bar{\boldsymbol{w}}$; and are defined as

$$\bar{\mathcal{V}} := \mathcal{V}; \; \bar{\mathcal{W}} := \mathcal{W} \tag{20}$$

where $\mathcal{V}$ and $\mathcal{W}$ are defined earlier in equation (5) and equation (6) respectively. Let $\mathcal{V}'$ be the function space for both the fine-scale component of the velocity $\boldsymbol{v}'$ and its corresponding weighting function $\boldsymbol{w}'$, and is defined as

$$\mathcal{V}' := \{\boldsymbol{v} \mid \boldsymbol{v} \in (H^1(\Omega^e))^{nd}, \, \boldsymbol{v} = \boldsymbol{0} \text{ on } \Gamma^e, e = 1, \ldots, N\} \tag{21}$$

The velocity field $\boldsymbol{v}(\boldsymbol{x})$ is now an element of the function space generated by the direct sum of $\bar{\mathcal{V}}$ and $\mathcal{V}'$, denoted by $\bar{\mathcal{V}} \oplus \mathcal{V}'$. Similarly the direct sum of $\bar{\mathcal{W}}$ and $\mathcal{V}'$, denoted by $\bar{\mathcal{W}} \oplus \mathcal{V}'$, is the function space for the field $\boldsymbol{w}(\boldsymbol{x})$.

In theory, we could decompose the pressure field into coarse-scale and fine-scale components. However, for simplicity we assume that there are no fine-scale terms for the pressure $p(\boldsymbol{x})$ and for its corresponding weighting function $q(\boldsymbol{x})$. Hence, the function space for the fields $p(\boldsymbol{x})$ and $q(\boldsymbol{x})$ is $\mathcal{P}$.

3.2. **Two-level classical mixed formulation.** The substitution of equations (17) and (18) into the classical mixed formulation, given by equations (12) and (13), becomes the first point of departure from the classical Galerkin formulation.

$$(\bar{w} + w', (\bar{v} + v') \cdot \nabla(\bar{v} + v')) + (\nabla(\bar{w} + w'), 2\nu\nabla(\bar{v} + v'))$$
$$-(\nabla \cdot (\bar{w} + w'), p) = (\bar{w} + w', b) + (\bar{w} + w', h)_{\Gamma_t} \quad (22)$$

$$(q, \nabla \cdot (\bar{v} + v')) = 0 \quad (23)$$

Because the weighting functions $\bar{w}$ and $w'$ are arbitrary, and because the functionals are linear in the weighting functions, we can write the above problem as two sub-problems. The *coarse-scale problem* can be written as:

$$(\bar{w}, (\bar{v} + v') \cdot \nabla(\bar{v} + v')) + (\nabla\bar{w}, 2\nu\nabla(\bar{v} + v')) - (\nabla \cdot \bar{w}, p)$$
$$= (\bar{w}, b) + (\bar{w}, h)_{\Gamma_t} \quad \forall \, \bar{w} \in \bar{\mathcal{W}} \quad (24)$$

$$(q, \nabla \cdot (\bar{v} + v')) = 0 \quad \forall \, q \in \mathcal{P} \quad (25)$$

The *fine-scale problem* can be written as:

$$(w', (\bar{v} + v') \cdot \nabla(\bar{v} + v')) + (\nabla w', 2\nu\nabla(\bar{v} + v')) - (\nabla \cdot w', p)$$
$$= (w', b) + (w', h)_{\Gamma_t} \quad \forall \, w' \in \mathcal{W}' \quad (26)$$

**Remark 2.** *Note that the fine scale problem is independent and uncoupled at the element level (defined over the sum of element interiors). Due to the assumption that the subgrid scale response vanishes on the element boundaries, $(\bar{w}, \bar{v} + v') = (\bar{w}, \bar{v}) + (\bar{w}, v')$.*

3.3. **Fine–scale interpolation and bubble functions.** If one chooses a single bubble function for interpolating the fine-scale variables (similar to the MINI element [12]), then we have

$$v' = b^e \beta; \quad w' = b^e \gamma \quad (27)$$

where $b^e$ is a bubble function, and $\beta$ and $\gamma$ are constant vectors. The gradients of the fine-scale velocity and weighting functions are

$$\nabla v' = \beta \nabla b^{eT}; \quad \nabla w' = \gamma \nabla b^{eT} \quad (28)$$

where $\nabla b^e$ is a dim $\times$ 1 vector of the derivatives of the bubble function. Standard bubble functions for several elements are provided in Table 1.

Table 1. Bubble Functions for Standard Finite Elements

| Element | Bubble function |
|---------|-----------------|
| T3      | $\xi_1\xi_2(1-\xi_1-\xi_2)$ |
| TET4    | $\xi_1\xi_2\xi_3(1-\xi_1-\xi_2-\xi_3)$ |
| Q4      | $(1-\xi_1^2)(1-\xi_2^2)$ |
| B8      | $(1-\xi_1^2)(1-\xi_2^2)(1-\xi_3^2)$ |

## 4. FIXED–POINT ITERATION

Equation (1) is nonlinear as engendered by the convection term $\boldsymbol{v}\cdot\nabla\boldsymbol{v}$. One popular method for dealing with the nonlinearity of this term is to linearize the weak form using fixed-point iteration. For a derivation of the linearization using fixed-point iteration, see Appendix 8.2. The linearized forms of the *fine* and *coarse-scale subproblems* are written as:

*Linearized coarse-scale subproblem*

$$(\bar{\boldsymbol{w}}, \boldsymbol{v}^c\cdot\nabla(\bar{\boldsymbol{v}}+\boldsymbol{v}')) + (\bar{\boldsymbol{w}}, (\bar{\boldsymbol{v}}+\boldsymbol{v}')\cdot\nabla\boldsymbol{v}^c) + (\nabla\bar{\boldsymbol{w}}, 2\nu\nabla(\bar{\boldsymbol{v}}+\boldsymbol{v}')) - (\nabla\cdot\bar{\boldsymbol{w}}, p)$$
$$= (\bar{\boldsymbol{w}}, \boldsymbol{b}) + (\bar{\boldsymbol{w}}, \boldsymbol{h})_{\Gamma_t} \quad (29)$$

$$(q, \nabla\cdot(\bar{\boldsymbol{v}}+\boldsymbol{v}')) = 0 \quad (30)$$

*Linearized fine-scale subproblem*

$$(\boldsymbol{w}', \boldsymbol{v}^c\cdot\nabla(\bar{\boldsymbol{v}}+\boldsymbol{v}')) + (\boldsymbol{w}', (\bar{\boldsymbol{v}}+\boldsymbol{v}')\cdot\nabla\boldsymbol{v}^c) + (\nabla\boldsymbol{w}', 2\nu\nabla(\bar{\boldsymbol{v}}+\boldsymbol{v}')) - (\nabla\cdot\boldsymbol{w}', p)$$
$$= (\boldsymbol{w}', \boldsymbol{b}) \quad (31)$$

To derive our stabilization parameter, we eliminate the fine-scale variables by solving the *fine-scale problem* (equation (31)) in terms of the coarse-scale variables. We then substitute the fine-scale solution into the *coarse-scale problem* (equation (29)) and solve the *coarse-scale problem* to obtain $\bar{\boldsymbol{v}}(\boldsymbol{x})$ and $p(\boldsymbol{x})$. Typically, the stabilization parameter is derived in a consistent manner by incorporating the coarse-scale residual evaluated over the element. Examples of such formulations include the work of Masud and Khurram [13] for the Navier–Stokes equations and that of Nakshatrala *et al* [14] for nearly incompressible linear elasticity. The derivation, which is presented in Appendix

8.3, leads to a stabilized weak formulation of the following form:

$$(\boldsymbol{w}, \boldsymbol{v}^c \cdot \nabla \boldsymbol{v}) + (\boldsymbol{w}, \boldsymbol{v} \cdot \nabla \boldsymbol{v}^c) + (\nabla \boldsymbol{w}, 2\nu \nabla \boldsymbol{v}) - (\nabla \cdot \boldsymbol{w}, p) + (q, \nabla \cdot \boldsymbol{v}) +$$
$$\left( \boldsymbol{v}^c \cdot \nabla \boldsymbol{w} + 2\nu \nabla^2 \boldsymbol{w} + \nabla q - (\nabla \boldsymbol{v}^c)^T \boldsymbol{w}, \boldsymbol{\tau}(\nabla p + \boldsymbol{v}^c \cdot \nabla \boldsymbol{v} + \boldsymbol{v} \cdot \nabla \boldsymbol{v}^c - 2\nu \nabla^2 \boldsymbol{v} - \boldsymbol{b}) \right)$$
$$= (\boldsymbol{w}, \boldsymbol{b}) + (\boldsymbol{w}, \boldsymbol{h})_{\Gamma_t} \quad (32)$$

**Remark 3.** *Notice that we have introduced the laplacian of the weighting function and the velocity as a result of integrating by parts. Computing these terms requires implementing the third order tensor which represents the second derivative of the shape functions.*

**Remark 4.** *Recall that these equations have been linearized about a certain velocity $\boldsymbol{v}^*$ with change in velocity $\delta \boldsymbol{v}$. For simplicity, we dropped the superscript on $\boldsymbol{v}^*$ and changed $\delta \boldsymbol{v}$ to $\boldsymbol{v}^c$, so to recover the stabilized nonlinear form of the weak formulation we must "reverse" the linearization process. Given that even the stabilization parameter $\boldsymbol{\tau}$ contains $\boldsymbol{v}^c$ and the complex expressions created from using integration by parts to isolate $\boldsymbol{v}'$, recovering the nonlinear form of Equation 32 remains difficult if not impossible.*

**Remark 5.** *Notice that the consistently linearized mixed variational multiscale form of the Navier Stokes equations is not self-adjoint as shown by the negative sign of the fourth term in the weighting slot of the stabilization term.*

The fixed-point iteration based weak form (equation (32)) is useful for low and moderate Reynolds flows, but the rate of convergence may be lower and the take a greater number of iterations [2]. In order to utilize homotopy methods (similar to arc-length), which greatly increase convergence for high Reynolds flows, one must posses the nonlinear weak formulation.

## 5. CONSISTENT NEWTON–RAPHSON SOLUTION STRATEGY

Instead of linearizing the convective term and solving the *coarse* and *fine scale subproblems* by substitution, as in the fixed point iteration strategy, for the consistent Newton–Raphson method, we treat equations (24)–(26) as global residuals and obtain the solution using an iterative update equation.

## 5.1. Scalar residual.
The scalar residual equations may be written as

$$r_c(\bar{v}; v', p) := (\bar{w}, (\bar{v} + v') \cdot \nabla(\bar{v} + v')) + (\nabla \bar{w}, 2\nu \nabla(\bar{v} + v'))$$
$$(33) \qquad - (\nabla \cdot \bar{w}, p) - (\bar{w}, b) - (\bar{w}, h)_{\Gamma_t}$$

$$(34) \qquad r_p(\bar{v}; v') := - (q, \nabla \cdot (\bar{v} + v'))$$

$$(35) \qquad r_f(\bar{v}; v', p) := (w', (\bar{v} + v') \cdot \nabla(\bar{v} + v')) + (\nabla w', 2\nu_f \nabla(\bar{v} + v')) - (\nabla \cdot w', p) - (w', b)$$

where the subscripts 'c', 'p', and 'f' stand for *coarse*, *pressure*, and *fine*.

## 5.2. Vector residual.
We again choose a single bubble function for interpolating the fine-scale variables. The solution, $\bar{v}$, and its weighting function, $\bar{w}$, may be expressed in terms of the nodal values $\hat{\bar{v}}$ and $\hat{\bar{w}}$ as

$$(36) \qquad \bar{v} = \hat{\bar{v}}^T N^T \; ; \; \bar{w} = \hat{\bar{w}}^T N^T$$

where $N$ is a row vector of shape functions for each node. Substituting equations (27) and (36) into equations (33)–(35) and noting the arbitrariness of $\hat{\bar{w}}$ and $\gamma$, we construct vector residuals, $R$, that are the sum contributions of the vector residuals at the element level, $R^e$, given as

$$R_c^e(\bar{v}; v', p) := \int_\Omega (N^T \odot I)((\bar{v} + v') \cdot \nabla)(\bar{v} + v') \, d\Omega$$
$$+ 2\nu \int_\Omega ((DN)J^{-1} \odot I)\text{vec}[\nabla \bar{v} + \nabla v'] \, d\Omega$$
$$(37) \qquad - \int_\Omega \text{vec}[J^{-T}(DN)^T]p \, d\Omega - \int_\Omega (N^T \odot I)b \, d\Omega - \int_{\Gamma_t} (N^T \odot I)h \, d\Gamma$$

$$(38) \qquad R_p^e(\bar{v}; v') := - \int_\Omega N^T \nabla \cdot (\bar{v} + v') \, d\Omega$$

$$R_f^e(\bar{v}; v', p) := \int_{\Omega_e} b^e I((\bar{v} + v') \cdot \nabla)(\bar{v} + v') \, d\Omega$$
$$+ 2\nu \int_{\Omega_e} (\nabla b^{eT} \odot I)\text{vec}[\nabla \bar{v} + \nabla v'] \, d\Omega$$
$$(39) \qquad - \int_{\Omega_e} (\nabla b^{eT} \odot I)\text{vec}[I]p \, d\Omega - \int_{\Omega_e} b^e Ib \, d\Omega$$

where $\boldsymbol{DN}$ represents a matrix of the first derivatives of the element shape functions, which is defined for a triangular element as

$$\boldsymbol{DN} := \begin{bmatrix} \frac{\partial N_1}{\partial \xi_1} & \frac{\partial N_1}{\partial \xi_2} \\ \vdots & \vdots \\ \frac{\partial N_3}{\partial \xi_1} & \frac{\partial N_3}{\partial \xi_2} \end{bmatrix} \tag{40}$$

$\boldsymbol{J}$ is the element jacobian matrix, vec[·] is an operation that represents a matrix with a vector, $\odot$ is the Kronecker product [15] (see Appendix 8.1).

5.3. **Extension to transient flows.** To extend the formulation presented above to transient flows, we need only to incorporate the time derivative of the velocity, $\dot{\boldsymbol{v}}$, in the balance of momentum, equation (1). This leads to the addition of the following terms, $(\bar{\boldsymbol{w}}, \dot{\boldsymbol{v}})$ and $(\boldsymbol{w}', \dot{\boldsymbol{v}})$ to the residual equations (33) and (35), respectively. This in turn adds the following terms, $\int_{\Omega}(\boldsymbol{N}^T \odot \boldsymbol{I})\dot{\boldsymbol{v}} \, d\Omega$ and $\int_{\Omega_e} b^e \boldsymbol{I}\dot{\boldsymbol{v}} \, d\Omega$ to the vector residuals, equations (37) and (39), respectively.

To integrate in time, from step $n$ to step $n+1$, we use the backward Euler method by which we express the acceleration, $\dot{\boldsymbol{v}}$, as the average of the velocity over time.

$$\dot{\boldsymbol{v}}_{n+1} = \frac{\boldsymbol{v}_{n+1} - \boldsymbol{v}_n}{\Delta t} \tag{41}$$

Substituting equation (41) into the additional acceleration terms in (37) and (39), we have

$$\int_{\Omega}(\boldsymbol{N}^T \odot \boldsymbol{I})\dot{\boldsymbol{v}}_{n+1} \, d\Omega = \int_{\Omega} \frac{1}{\Delta t}(\boldsymbol{N}^T \odot \boldsymbol{I})\bar{\boldsymbol{v}}_{n+1} \, d\Omega - \int_{\Omega} \frac{1}{\Delta t}(\boldsymbol{N}^T \odot \boldsymbol{I})\bar{\boldsymbol{v}}_n \, d\Omega \tag{42}$$

$$\int_{\Omega_e} b^e \boldsymbol{I}\dot{\boldsymbol{v}} \, d\Omega = \int_{\Omega_e} \frac{b^e}{\Delta t}\boldsymbol{I}\bar{\boldsymbol{v}}_{n+1} \, d\Omega - \int_{\Omega_e} \frac{b^e}{\Delta t}\boldsymbol{I}\bar{\boldsymbol{v}}_n \, d\Omega \tag{43}$$

where $\bar{\boldsymbol{v}}_n$ is the converged velocity from the previous time step.

5.4. **Tangent matrix.** Using a Newton-Raphson type approach, we obtain the solution in an iterative fashion using the following update equation until the residual is under a prescribed tolerance,

$$\bar{\boldsymbol{v}}^{i+1} = \bar{\boldsymbol{v}}^i + \Delta \bar{\boldsymbol{v}}^i \quad ; \quad p^{i+1} = p^i + \Delta p^i \quad ; \quad \boldsymbol{v}'^{i+1} = \boldsymbol{v}'^i + \Delta \boldsymbol{v}'^i \tag{44}$$

where the updates at each iteration, $i$ are calculated from the following system of equations

$$\begin{bmatrix} \frac{D\boldsymbol{R}_c}{D\bar{\boldsymbol{v}}} & \frac{D\boldsymbol{R}_c}{Dp} & \frac{D\boldsymbol{R}_c}{D\boldsymbol{v}'} \\ \frac{D\boldsymbol{R}_p}{D\bar{\boldsymbol{v}}} & \frac{D\boldsymbol{R}_p}{Dp} & \frac{D\boldsymbol{R}_p}{D\boldsymbol{v}'} \\ \frac{D\boldsymbol{R}_f}{D\bar{\boldsymbol{v}}} & \frac{D\boldsymbol{R}_f}{Dp} & \frac{D\boldsymbol{R}_f}{D\boldsymbol{v}'} \end{bmatrix} \begin{Bmatrix} \Delta \bar{\boldsymbol{v}} \\ \Delta p \\ \Delta \boldsymbol{v}' \end{Bmatrix} = - \begin{Bmatrix} \boldsymbol{R}_c \\ \boldsymbol{R}_p \\ \boldsymbol{R}_f \end{Bmatrix} \tag{45}$$

The element matrices, $\frac{\mathrm{D}\boldsymbol{R}^e}{\mathrm{D}(\cdot)}$, which are assembled to form the consistent tangent matrix, are derived in Appendix 8.4.

## 6. NUMERICAL RESULTS

In this section, we present the results from several numerical experiments that directly compare the Newton-Raphson solution approach to fixed-point iteration.

6.1. **Body force driven cavity.** We first perform a mesh refinement study to observe the convergence rates for the Newton-Raphson and fixed-point iteration strategies for a smooth problem using linear elements as the element length decreases. The problem consists of a body force driven cavity for which the solution may easily be obtained analytically. The geometry and a typical mesh, shown in Figure 1, consists of a unit domain with boundaries that prohibit flow in the tangential and normal directions ($v_x = v_y = 0.0$). The prescribed constant body force, $\boldsymbol{b}$, in components, is given as

$$b_x = (12 - 24y)x^4 + (-24 + 48y)x^3 + (-48y + 72y^2 - 48y^3 + 12)x^2$$
$$+ (-2 + 24y - 72y^2 + 48y^3)x + 1 - 4y + 12y^2 - 8y^3$$
$$+ (4xy - 12xy^2 + 8xy^3 - 12x^2y + 36x^2y^2 - 24x^2y^3 + 8x^3y - 24x^3y^2 + 16x^3y^3)*$$
$$(x^2(1-x)^2(2y - 6y^2 + 4y^3))$$
$$+ (2x^2 - 12x^2y + 12x^2y^2 - 4x^3 + 24x^3y - 24x^3y^2 + 2x^4 - 12x^4y + 12x^4y^2)*$$
$$(46) \qquad (-y^2(1-y)^2(2x - 6x^2 + 4x^3))$$
$$b_y = (8 - 48y + 48y^2)x^3 + (-12 + 72y - 72y^2)x^2$$
$$+ (4 - 24y + 48y^2 - 48y^3 + 24y^4)x - 12y^2 + 24y^3 - 12y^4$$
$$+ (-2y^2 + 12y^2x - 12y^2x^2 + 4y^3 - 24y^3x + 24y^3x^2 - 2y^4 + 12y^4x - 12y^4x^2)*$$
$$(x^2(1-x)^2(2y - 6y^2 + 4y^3))$$
$$+ (-4yx + 12yx^2 - 8yx^3 + 12y^2x - 36y^2x^2 + 24y^2x^3 - 8y^3x + 24y^3x^2 - 16y^3x^3)*$$
$$(47) \qquad (-y^2(1-y)^2(2x - 6x^2 + 4x^3))$$

For a unit viscosity, the exact solution is

$$v_x = x^2(1-x)^2(2y - 6y^2 + 4y^3)$$

$$v_y = -y^2(1-y)^2(2x - 6x^2 + 4x^3)$$

(48) $$p = x(1-x)$$

The results, for triangular elements of edge length $h$ with a prescribed body force and known solution, are shown in Figures 2 and 3. In Figure 2 the values of the pressure and velocity are shown as measured along the x-axis at $y = 0.5$. From Figure 3 we see that the $H^1$-semi norm of the pressure is 1.0 and the $L^2$ norm of the velocity is 2.0, respectively for both methods. We also see that the error in the velocity is slightly lower for the Newton-Raphson strategy.

In Figure 4 we show a comparison between the residual convergence for the consistent Newton-Schur approach and the fixed-point iteration method for the smooth body force problem at $Re = 400$. The solution converges in three steps for the Newton-Raphson approach, whereas after ten iterations, the fixed-point method has not yet reached the tolerance. For the higher Reynolds case, $Re = 5,000$, the Newton-Raphson strategy again takes three iterations to converge, but the fixed-point iteration method diverges. In this case, one would have to employ a continuation technique to obtain the solution using the fixed-point iteration method.

6.2. **Lid-driven cavity.** The lid-driven cavity has served as a benchmark problem for numerical methods and has been analyzed by a number of authors (for examples, see [16, 17, 18, 19]). The problem description, boundary conditions, and a typical mesh are shown in Figure 5. For a unit length and height domain, a unit velocity in the horizontal direction is prescribed across the top boundary, while the sides and bottom are treated as solid surfaces with no flow in the horizontal or vertical direction. To resolve the boundary conditions at the corners they are treated as part of the vertical walls. Therefor, the prescribed velocity approaches zero at the corners. Plots of the stream traces over the velocity contours are shown for $Re = 400$ and $Re = 5,000$ in Figures 7 and 8 respectively for the Newton-Raphson approach. Similar plots for the fixed-point iteration strategy have been presented in [13]. Results for both methods are compared to other published results in Figure 6. Figure 6 shows the velocity magnitude in the x-direction as measured along the y-axis at $x = 0.5$. The results closely match those presented in [16]. For the $Re = 5,000$ case, aside from the central vortex, the lower vortices and the upper left vortex are captured for a mesh of triangular elements with edge length 0.025 using the Newton-Raphson based approach. We were unable to obtain a converged solution to the lid-driven cavity problem using the fixed-point iteration scheme at $Re = 5,000$.

For the lid-driven cavity, we observe quadratic convergence as expected for the consistent Newton-Raphson approach. The residual is listed for each iteration in Table 2. Above $Re = 5,000$, we were unable to obtain a converged solution for the Newton-Raphson based approach. Even when a continuation method is employed by solving the problem at low $Re$ and slowly increasing $Re$ by a factor of 1.1, using the previous converged solution as the starting point, the method does not converge. There seems to exist a well defined $Re$ above which the consistent method does not converge. The lack of convergence for consistent methods at very high $Re$ is pointed out in [2].

TABLE 2. Quadratic Convergence for Lid–Driven Cavity Problem at $Re = 400$

| Iteration | Residual |
|---|---|
| 1 | 1.429492E+04 |
| 2 | 8.608793E+03 |
| 3 | 4.716364E+03 |
| 4 | 2.104800E+03 |
| 5 | 6.737382E+02 |
| 6 | 5.754620E+01 |
| 7 | 4.170450E-01 |
| 8 | 5.516418E-04 |
| 9 | 4.784608E-08 |
| 10 | 1.784712E-11 |

6.3. **Backward facing step.** Another well-known benchmark problem that posses a corner singularity is the backward facing step. The domain, boundary conditions, and mesh are shown in Figure 9. A parabolic velocity is prescribed over the inflow face with a maximum value of 1.0 at the center of the opening. The results for the consistent Newton-Raphson method are presented in Figure 10, which show streamtraces over the velocity contours and the pressure contours respectively for $Re = 150$. Results for the fixed-point iteration method are presented in [13]. In order to obtain the solution for both the fixed-point iteration and consistent Newton-Raphson approaches, a continuation process must be employed. To compare the convergence of the residual for both cases, we present the results for the first converged Reynolds number in the continuation. Figure 11 shows the convergence of the residual at $Re = 15$. Again, the Newton-Raphson method converges in much fewer iterations.

### 6.4. Consistent Newton-Raphson and the transition from steady to periodic flow.

Whereas the convergence for the Newton-Raphson based approach is much faster than the fixed-point iteration based method, the Newton-Raphson based approach exhibits vulnerability in the region of turning or bifurcation points. In particular, the consistent method does not converge when the flow past a bluff object transitions from steady to periodic vortex shedding. This transition represents an instance of Hopf bifurcation [20].

The domain and boundary conditions for two problems that engender this behavior are shown in Figures 12 and 14 for flow past a circular and square cylinder respectively. The computational meshes used for this problem are also shown in Figures 12 and 14. Figures 13 and 15 show the last converged iterations for the consistent Newton-Raphson based approach immediately before the method fails to converge. Notice that for both cases, periodic vortex shedding has not yet started. The residual for each case is shown in Tables 3 and 4.

The consistent Newton-Raphson method does not fail for all transient vortex shedding related problems, for example, the jet flow through an orifice problem described in Figure 16. The computational mesh is also shown in Figure 16. Figure 17 shows the results using the consistent Newton-Raphson approach for a jet flowing through a 1/16 unit opening in an infinite wall for viscosity, $\nu = 0.005$, and maximum flow speed at the orifice of $\boldsymbol{v} = (1.0, 0)$. The results closely match those reported in [5] for the same Reynolds number. Again quadratic convergence is obtained for the Newton-Raphson technique for an initial guess of $\boldsymbol{v} = 0$. For the jet flow problem, there is no transition from steady to periodic flow.

### 7. CONCLUSIONS

We have presented a comparison between a fixed-point iteration and consistent Newton-Raphson solution strategy for a variational multiscale formulations for incompressible Navier–Stokes. The Newton-Raphson approach shows faster convergence for several numerical test problems and in many cases converges quadratically in the vicinity of the solution. We have shown that in some instances, the Newton-Raphson approach converges while the fixed-point iteration strategy does not. We have also shown that for some problems, the Newton-Raphson strategy is unable to pass through the transition from steady to periodic flow without additional sophistication being built into the algorithm. In summary, this work points to the need for advanced solution algorithms able to maintain the convergence qualities of the Newton-Raphson technique with increased robustness for high Reynolds flows and flows with distinct bifurcation points.

TABLE 3. Convergence for Flow Past a Circular Cylinder at $Re = 75$

| Circular Cylinder | | |
| --- | --- | --- |
| | Prior to Hopf Bifurcation | At Hopf Bifurcation |
| Iteration | Residual | Residual |
| 1 | 7.338022E+01 | 6.066543E+01 |
| 2 | 5.783447E-02 | 2.542188E-01 |
| 3 | 4.012763E-02 | 2.447873E-01 |
| 4 | 4.246893E-07 | 2.451840E-01 |
| 5 | | 2.451840E-01 |
| 6 | | 2.451840E-01 |
| 7 | | . . . |

TABLE 4. Convergence for Flow Past a Square Cylinder at $Re = 75$

| Square Cylinder | | |
| --- | --- | --- |
| | Prior to Hopf Bifurcation | At Hopf Bifurcation |
| Iteration | Residual | Residual |
| 1 | 2.526981E+01 | 1.460652E+00 |
| 2 | 1.561594E-02 | 7.228785E-02 |
| 3 | 5.892473E-03 | 7.196581E-02 |
| 4 | 4.721279E-08 | 7.195860E-02 |
| 5 | | 7.195860E-02 |
| 6 | | 7.195860E-02 |
| 7 | | . . . |

## 8. APPENDIX

## 8.1. Notation and definitions.
Consider an $n \times m$ matrix $\boldsymbol{A}$ and a $p \times q$ matrix $\boldsymbol{B}$

$$\boldsymbol{A} = \begin{bmatrix} a_{1,1} & \cdots & a_{1,m} \\ \vdots & \ddots & \vdots \\ a_{n,1} & \cdots & a_{n,m} \end{bmatrix} ; \boldsymbol{B} = \begin{bmatrix} b_{1,1} & \cdots & b_{1,q} \\ \vdots & \ddots & \vdots \\ b_{p,1} & \cdots & b_{p,q} \end{bmatrix}$$

The *Kronecker product* of these matrices is an $np \times mq$ matrix, and is defined as

$$\boldsymbol{A} \odot \boldsymbol{B} := \begin{bmatrix} a_{1,1}\boldsymbol{B} & \cdots & a_{1,m}\boldsymbol{B} \\ \vdots & \ddots & \vdots \\ a_{n,1}\boldsymbol{B} & \cdots & a_{n,m}\boldsymbol{B} \end{bmatrix}$$

The vec[·] operator is defined as

$$\mathrm{vec}[\boldsymbol{A}] := \begin{bmatrix} a_{1,1} \\ \vdots \\ a_{1,m} \\ \vdots \\ a_{n,1} \\ \vdots \\ a_{n,m} \end{bmatrix}$$

## 8.2. Fixed–point linearization of the convective term.
We linearize the convective term, $\boldsymbol{v} \cdot \nabla \boldsymbol{v}$ in the weak form about a velocity $\boldsymbol{v}^*$ with a change in velocity $\delta \boldsymbol{v}$ as follows.

$$\boldsymbol{r}(\boldsymbol{v}) = (\boldsymbol{w}, \boldsymbol{v} \cdot \nabla \boldsymbol{v})$$

$$\frac{\mathrm{d}}{\mathrm{d}\epsilon}\Big[\boldsymbol{r}(\boldsymbol{v}^* + \epsilon \delta \boldsymbol{v})\Big]_{\epsilon=0} = \frac{\mathrm{d}}{\mathrm{d}\epsilon}\Big[(\boldsymbol{w}, (\boldsymbol{v}^* + \epsilon \delta \boldsymbol{v}) \cdot \nabla (\boldsymbol{v}^* + \epsilon \delta \boldsymbol{v}))\Big]_{\epsilon=0}$$

$$= \Big[(\boldsymbol{w}, \delta \boldsymbol{v} \cdot \nabla (\boldsymbol{v}^* + \epsilon \delta \boldsymbol{v})) + (\boldsymbol{w}, (\boldsymbol{v}^* + \epsilon \delta \boldsymbol{v}) \cdot \nabla \delta \boldsymbol{v})\Big]_{\epsilon=0}$$

(49)
$$= (\boldsymbol{w}, \delta \boldsymbol{v} \cdot \nabla \boldsymbol{v}^*) + (\boldsymbol{w}, \boldsymbol{v}^* \cdot \nabla \delta \boldsymbol{v})$$

To simplify the formulation, we drop the superscript on $\boldsymbol{v}^*$ and refer to the change in velocity $\delta \boldsymbol{v}$ as $\boldsymbol{v}^c$. The linearized form of the second term on the left side of equation (12) is now expressed as

(50)
$$(\boldsymbol{w}, \boldsymbol{v} \cdot \nabla \boldsymbol{v}) = (\boldsymbol{w}, \boldsymbol{v}^c \cdot \nabla \boldsymbol{v}) + (\boldsymbol{w}, \boldsymbol{v} \cdot \nabla \boldsymbol{v}^c)$$

**8.3. Derivation of the stabilization term for the fixed-point iteration strategy.** Due to the linearity of the solution and using integration by parts on the pressure term, we re-arrange equation (31) as

$$(w', v^c \cdot \nabla v') + (w', v' \cdot \nabla v^c) + (\nabla w', 2\nu \nabla v') = (w', b_f) - (w', \nabla p)$$
(51)
$$-(w', \dot{\bar{v}}) - (w', v^c \cdot \nabla \bar{v}) - (w', \bar{v} \cdot \nabla v^c) - (\nabla w', 2\nu \nabla \bar{v})$$

Performing integration by parts on the last term in the previous equation and collecting terms we reduce this expression to

(52)
$$(w', v^c \cdot \nabla v') + (w', v' \cdot \nabla v^c) + (\nabla w', 2\nu \nabla v') = (w', b_f - \nabla p - \dot{\bar{v}} - v^c \cdot \nabla \bar{v} - \bar{v} \cdot \nabla v^c + 2\nu \nabla^2 \bar{v})$$

Making the substitution $\bar{r} = b_f - \nabla p - \dot{\bar{v}} - v^c \cdot \nabla \bar{v} - \bar{v} \cdot \nabla v^c + 2\nu \nabla^2 \bar{v}$, we have

(53)
$$(w', v^c \cdot \nabla v') + (w', v' \cdot \nabla v^c) + (\nabla w', 2\nu \nabla v') = (w', \bar{r})$$

We now solve equation (53) analytically using bubble functions. We assume that $v'$ and $w'$ are represented by bubble functions over $\Omega'_f$, as in equation (27), and move the constant coefficients outside of the integrals to obtain

(54)
$$\gamma^T \left( \int_{\Omega^e_f} (b^e v^c \cdot \nabla b^e + \nu |\nabla b^e|^2) \mathbf{I} + (b^e)^2 \nabla v^c + \nu \nabla b^e \otimes \nabla b^e \, d\Omega \right) \beta = \gamma^T (b^e, \bar{r})$$

where $\nabla v^c$ is the gradient of the known velocity about which the equations have been linearized and $\nabla b^e$ is a dim $\times$ 1 vector of the derivatives of the bubble functions. Because the coefficients $\gamma$ are arbitrary, we have a system of equations that we can solve for $\beta$ in terms of $\bar{r}$,

(55)
$$\beta = \mathbf{A}^{-1} \mathbf{R}$$

where $\mathbf{A}$ and $\mathbf{R}$ are defined as follows

(56)
$$\mathbf{A} = \int_{\Omega^e_f} (b^e v^c \cdot \nabla b^e + \nu |\nabla b^e|^2) \mathbf{I} + (b^e)^2 \nabla v^c + \nu \nabla b^e \otimes \nabla b^e \, d\Omega$$

(57)
$$\mathbf{R} = \int_{\Omega^e_f} b^e \bar{r} \, d\Omega$$

we now have an expression for the fine-scale velocity which we can substitute into the *coarse-scale subproblem*.

(58)
$$v'(x) = b^e \beta = b^e \mathbf{A}^{-1} \mathbf{R}$$

In the finite element setting, as the mesh becomes adequately refined, $\bar{r}$ will essentially be constant over an element. We can then move $\bar{r}$ to the outside of the integral in equation (57) and rewrite our approximation for the fine-scale velocity as

$$v'(x) = b^e \int_{\Omega_f^e} b^e \, d\Omega \mathbf{A}^{-1} \bar{r} = \tau \bar{r} \tag{59}$$

where the stabilization parameter $\tau$ (in this case a second-order tensor) is defined as

$$\tau = b^e \int_{\Omega_f^e} b^e \, d\Omega \mathbf{A}^{-1}$$

$$= b^e \int_{\Omega_f^e} b^e \, d\Omega \left[ \int_{\Omega_f^e} (b^e v^c \cdot \nabla b^e + \nu |\nabla b^e|^2) \mathbf{I} + (b^e)^2 \nabla v^c + \nu \nabla b^e \otimes \nabla b^e \, d\Omega \right]^{-1} \tag{60}$$

Using integration by parts, the linearity of the solution field, and the identity, $(\mathbf{a}, \mathbf{A}\mathbf{b}) = (\mathbf{A}^T \mathbf{a}, \mathbf{b})$, we substitute equation (59) into the *coarse scale subproblem* to get

$$(\boldsymbol{w}, \dot{\boldsymbol{v}}) + (\boldsymbol{w}, \boldsymbol{v}^c \cdot \nabla \boldsymbol{v}) + (\boldsymbol{w}, \boldsymbol{v} \cdot \nabla \boldsymbol{v}^c) + (\nabla \boldsymbol{w}, 2\nu \nabla \boldsymbol{v}) - (\nabla \cdot \boldsymbol{w}, p) + (q, \nabla \cdot \boldsymbol{v}) +$$

$$\left( \boldsymbol{v}^c \cdot \nabla \boldsymbol{w} + 2\nu \nabla^2 \boldsymbol{w} + \nabla q - (\nabla \boldsymbol{v}^c)^T \boldsymbol{w}, \tau(\nabla p + \dot{\boldsymbol{v}} + \boldsymbol{v}^c \cdot \nabla \boldsymbol{v} + \boldsymbol{v} \cdot \nabla \boldsymbol{v}^c - 2\nu \nabla^2 \boldsymbol{v} - \boldsymbol{b}) \right)$$

$$= (\boldsymbol{w}, \boldsymbol{b}) + (\boldsymbol{w}, \boldsymbol{h})_{\Gamma_t} \tag{61}$$

**8.4. Derivation of the terms for the consistent tangent matrix.**

$$\frac{D\boldsymbol{R}_c}{D\bar{\boldsymbol{v}}} \cdot \delta\bar{\boldsymbol{v}} = \int_{\Omega} (\boldsymbol{N}^T \odot \boldsymbol{I}) \delta\dot{\bar{\boldsymbol{v}}} \, d\Omega + \int_{\Omega} (\boldsymbol{N}^T \odot \boldsymbol{I})((\delta\bar{\boldsymbol{v}} \cdot \nabla)(\bar{\boldsymbol{v}} + \boldsymbol{v}') + (\bar{\boldsymbol{v}} + \boldsymbol{v}') \cdot \nabla \delta\bar{\boldsymbol{v}}) \, d\Omega$$

$$+ 2\nu_f \int_{\Omega} ((\boldsymbol{DN})\boldsymbol{J}^{-1} \odot \boldsymbol{I}) \text{vec}[\nabla \delta\bar{\boldsymbol{v}}] \, d\Omega \tag{62}$$

$$\frac{D\boldsymbol{R}_p}{D\bar{\boldsymbol{v}}} \cdot \delta\bar{\boldsymbol{v}} = -\int_{\Omega} \boldsymbol{N}^T \nabla \cdot \delta\bar{\boldsymbol{v}} \, d\Omega \tag{63}$$

$$\frac{D\boldsymbol{R}_c}{D\boldsymbol{v}'} \cdot \delta\boldsymbol{v}' = \int_{\Omega} (\boldsymbol{N}^T \odot \boldsymbol{I})((\delta\boldsymbol{v}' \cdot \nabla)(\bar{\boldsymbol{v}} + \boldsymbol{v}') + (\bar{\boldsymbol{v}} + \boldsymbol{v}') \cdot \nabla \delta\boldsymbol{v}') \, d\Omega$$

$$+ 2\nu_f \int_{\Omega} ((\boldsymbol{DN})\boldsymbol{J}^{-1} \odot \boldsymbol{I}) \text{vec}[\nabla \delta\boldsymbol{v}'] \, d\Omega \tag{64}$$

$$\frac{D\boldsymbol{R}_p}{D\boldsymbol{v}'} \cdot \delta\boldsymbol{v}' = -\int_{\Omega} \boldsymbol{N}^T \nabla \cdot \delta\boldsymbol{v}' \, d\Omega \tag{65}$$

$$\frac{D\boldsymbol{R}_c}{Dp} \cdot \delta p = -\int_{\Omega} \text{vec}[\boldsymbol{J}^{-T}(\boldsymbol{DN})^T] \delta p \, d\Omega \tag{66}$$

$$\frac{\mathrm{D}\bm{R}_f}{\mathrm{D}\bar{\bm{v}}} \cdot \delta\bar{\bm{v}} = \int_{\Omega_e} b^e \bm{I}\delta\dot{\bar{\bm{v}}} \, \mathrm{d}\Omega + \int_{\Omega_e} b^e \bm{I}(\delta\bar{\bm{v}} \cdot \nabla)(\bar{\bm{v}} + \bm{v}') + ((\bar{\bm{v}} + \bm{v}') \cdot \nabla)\delta\bar{\bm{v}} \, \mathrm{d}\Omega$$

$$(67) \qquad + \, 2\nu_f \int_{\Omega_e} (\nabla b^{eT} \odot \bm{I}) \mathrm{vec}[\nabla \delta\bar{\bm{v}}] \, \mathrm{d}\Omega$$

$$\frac{\mathrm{D}\bm{R}_f}{\mathrm{D}\bm{v}'} \cdot \delta\bm{v}' = \int_{\Omega_e} b^e \bm{I}(\delta\bm{v}' \cdot \nabla)(\bar{\bm{v}} + \bm{v}') + ((\bar{\bm{v}} + \bm{v}') \cdot \nabla)\delta\bm{v}' \, \mathrm{d}\Omega$$

$$(68) \qquad + \, 2\nu_f \int_{\Omega_e} (\nabla b^{eT} \odot \bm{I}) \mathrm{vec}[\nabla \delta\bm{v}'] \, \mathrm{d}\Omega$$

$$(69) \qquad \frac{\mathrm{D}\bm{R}_f}{\mathrm{D}p} \cdot \delta p = -\int_{\Omega_e} (\nabla b^{eT} \odot \bm{I}) \mathrm{vec}[\bm{I}] \delta p$$

In the finite element setting, we introduce a finite element discretization of the solution for the pressure and velocity as follows

$$(70) \qquad \bar{\bm{v}} = \hat{\bar{\bm{v}}}^T \bm{N}^T \quad ; \quad \bm{v}' = b^e(\xi)\bm{\beta} \quad ; \quad p = \bm{N}\hat{\bm{p}}$$

Making use of the following equivalences:

$$(71) \qquad \frac{\mathrm{D}(\cdot)}{\mathrm{D}\bar{\bm{v}}} \cdot \delta\bar{\bm{v}} \equiv \frac{\mathrm{D}(\cdot)}{\mathrm{D}\hat{\bar{\bm{v}}}} \cdot \mathrm{vec}[\delta\hat{\bar{\bm{v}}}^T] \quad ; \quad \frac{\mathrm{D}(\cdot)}{\mathrm{D}\bm{v}'} \cdot \delta\bm{v}' \equiv \frac{\mathrm{D}(\cdot)}{\mathrm{D}\bm{\beta}} \cdot \delta\bm{\beta} \quad ; \quad \frac{\mathrm{D}(\cdot)}{\mathrm{D}p} \cdot \delta p \equiv \frac{\mathrm{D}(\cdot)}{\mathrm{D}\hat{\bm{p}}} \cdot \delta\hat{\bm{p}}$$

we derive the the discretized form of the tangent stiffness matrix as follows,

$$\frac{\mathrm{D}\bm{R}_c^e}{\mathrm{D}\hat{\bar{\bm{v}}}} \cdot \mathrm{vec}[\delta\hat{\bar{\bm{v}}}^T] = \int_{\Omega} (\bm{N}^T\bm{N} \odot \bm{I})\frac{1}{\Delta t}\mathrm{vec}[\delta\hat{\bar{\bm{v}}}^T] \, \mathrm{d}\Omega + \int_{\Omega} (\bm{N}^T \odot \bm{I})\nabla(\bar{\bm{v}} + \bm{v}')(\bm{N} \odot \bm{I})\mathrm{vec}[\delta\hat{\bar{\bm{v}}}^T] \, \mathrm{d}\Omega$$

$$+ \int_{\Omega} (\bm{N}^T(\bar{\bm{v}} + \bm{v}')^T \bm{J}^{-T}(\bm{D}\bm{N})^T \odot \bm{I})\mathrm{vec}[\delta\hat{\bar{\bm{v}}}^T] \, \mathrm{d}\Omega$$

$$(72) \qquad + \, 2\nu_f \int_{\Omega} ((\bm{D}\bm{N})\bm{J}^{-1}\bm{J}^{-T}(\bm{D}\bm{N})^T \odot \bm{I})\mathrm{vec}[\delta\hat{\bar{\bm{v}}}^T] \, \mathrm{d}\Omega$$

$$(73) \qquad \frac{\mathrm{D}\bm{R}_p^e}{\mathrm{D}\hat{\bar{\bm{v}}}} \cdot \mathrm{vec}[\delta\hat{\bar{\bm{v}}}^T] = -\int_{\Omega} \bm{N}^T \mathrm{vec}(\bm{J}^{-T}(\bm{D}\bm{N})^T)^T \mathrm{vec}[\delta\hat{\bar{\bm{v}}}^T] \, \mathrm{d}\Omega$$

$$\frac{\mathrm{D}\bm{R}_c^e}{\mathrm{D}\bm{\beta}} \cdot \delta\bm{\beta} = \int_{\Omega} (\bm{N}^T \odot \bm{I})b^e \nabla(\bar{\bm{v}} + \bm{v}') \, \delta\bm{\beta} \, \mathrm{d}\Omega + \int_{\Omega} (\bm{N}^T \odot \bm{I})\nabla b^{eT}(\bar{\bm{v}} + \bm{v}')\bm{I} \, \delta\bm{\beta} \, \mathrm{d}\Omega$$

$$(74) \qquad + \, 2\nu_f \int_{\Omega} ((\bm{D}\bm{N})\bm{J}^{-1}\nabla b^e \odot \bm{I}) \, \delta\bm{\beta} \, \mathrm{d}\Omega$$

$$(75) \qquad \frac{\mathrm{D}\bm{R}_p^e}{\mathrm{D}\bm{\beta}} \cdot \delta\bm{\beta} = -\int_{\Omega} \bm{N}^T \nabla b^{eT} \, \delta\bm{\beta} \, \mathrm{d}\Omega$$

$$(76) \qquad \frac{\mathrm{D}\bm{R}_c^e}{\mathrm{D}\hat{\bm{p}}} \cdot \delta\hat{\bm{p}} = -\int_{\Omega} \mathrm{vec}[\bm{J}^{-T}(\bm{D}\bm{N})^T]\bm{N}\delta\hat{\bm{p}} \, \mathrm{d}\Omega$$

$$(77) \quad \frac{\mathrm{D}\boldsymbol{R}_f^e}{\mathrm{D}\hat{\bar{\boldsymbol{v}}}} \cdot \mathrm{vec}[\delta\hat{\bar{\boldsymbol{v}}}^T] = \int_{\Omega_e} b^e \boldsymbol{I}(\boldsymbol{N} \odot \boldsymbol{I})\frac{1}{\Delta t}\mathrm{vec}[\delta\hat{\bar{\boldsymbol{v}}}^T]\,\mathrm{d}\Omega + \int_{\Omega_e} b^e \boldsymbol{I}\nabla(\bar{\boldsymbol{v}} + \boldsymbol{v}')(\boldsymbol{N} \odot \boldsymbol{I})\mathrm{vec}[\delta\hat{\bar{\boldsymbol{v}}}^T]\,\mathrm{d}\Omega$$
$$+ \int_{\Omega_e} b^e \boldsymbol{I}((\bar{\boldsymbol{v}} + \boldsymbol{v}')^T \boldsymbol{J}^{-T}(\boldsymbol{DN})^T \odot \boldsymbol{I})\mathrm{vec}[\delta\hat{\bar{\boldsymbol{v}}}^T]\,\mathrm{d}\Omega$$
$$+ 2\nu_f \int_{\Omega_e} (\nabla b^{eT} \boldsymbol{J}^{-T}(\boldsymbol{DN})^T \odot \boldsymbol{I})\mathrm{vec}[\nabla\delta\bar{\boldsymbol{v}}]\,\mathrm{d}\Omega$$

$$(78) \quad \frac{\mathrm{D}\boldsymbol{R}_f^e}{\mathrm{D}\boldsymbol{\beta}} \cdot \delta\boldsymbol{\beta} = \int_{\Omega_e} b^e \boldsymbol{I}(b^e \nabla(\bar{\boldsymbol{v}} + \boldsymbol{v}') + (\nabla b^e \cdot (\bar{\boldsymbol{v}} + \boldsymbol{v}'))\boldsymbol{I})\delta\boldsymbol{\beta}\,\mathrm{d}\Omega + 2\nu_f \int_{\Omega_e} (\nabla b^{eT}\nabla b^e \odot \boldsymbol{I})\delta\boldsymbol{\beta}\,\mathrm{d}\Omega$$

$$(79) \quad \frac{\mathrm{D}\boldsymbol{R}_f^e}{\mathrm{D}\hat{p}} \cdot \delta\hat{p} = -\int_{\Omega_e} (\nabla b^{eT} \odot \boldsymbol{I})\mathrm{vec}[\boldsymbol{I}]\boldsymbol{N}\delta\hat{p}\,\mathrm{d}\Omega$$

Due to the nontrivial nature of the derivation of $\frac{\mathrm{D}\boldsymbol{R}_f^e}{\mathrm{D}\boldsymbol{\beta}} \cdot \delta\boldsymbol{\beta}$, the derivation is presented as follows:

$$\frac{\mathrm{D}\boldsymbol{R}_e'}{\mathrm{D}\boldsymbol{\beta}} \cdot \delta\boldsymbol{\beta} = \int_{\Omega_e} b^e \boldsymbol{I}(b^e \delta\boldsymbol{\beta} \cdot \nabla)(\bar{\boldsymbol{v}} + \boldsymbol{v}') + ((\bar{\boldsymbol{v}} + \boldsymbol{v}') \cdot \nabla)(b^e \delta\boldsymbol{\beta})\,\mathrm{d}\Omega$$
$$+ 2\nu_f \int_{\Omega_e} (\nabla b^{eT} \odot \boldsymbol{I})(\nabla b^e \odot \boldsymbol{I})\delta\boldsymbol{\beta}\,\mathrm{d}\Omega$$
$$= \int_{\Omega_e} b^e \boldsymbol{I}(b^e \delta\boldsymbol{\beta} \cdot \nabla)(\bar{\boldsymbol{v}} + \boldsymbol{v}') + (\delta\boldsymbol{\beta} \otimes \nabla b^e) \cdot (\bar{\boldsymbol{v}} + \boldsymbol{v}')\,\mathrm{d}\Omega$$
$$+ 2\nu_f \int_{\Omega_e} (\nabla b^{eT} \odot \boldsymbol{I})(\nabla b^e \odot \boldsymbol{I})\delta\boldsymbol{\beta}\,\mathrm{d}\Omega$$
$$= \int_{\Omega_e} b^e \boldsymbol{I}(b^e \nabla(\bar{\boldsymbol{v}} + \boldsymbol{v}') + (\nabla b^e \cdot (\bar{\boldsymbol{v}} + \boldsymbol{v}'))\boldsymbol{I})\delta\boldsymbol{\beta}\,\mathrm{d}\Omega$$
$$(80) \quad + 2\nu_f \int_{\Omega_e} (\nabla b^{eT}\nabla b^e \odot \boldsymbol{I})\delta\boldsymbol{\beta}\,\mathrm{d}\Omega$$

Noting the arbitrariness of $\mathrm{vec}[\delta\hat{\bar{\boldsymbol{v}}}^T]$, $\delta\boldsymbol{\beta}$, and $\delta\hat{p}$, we have expressions for the tangent stiffness matrix components in the finite element setting:

$$\frac{\mathrm{D}\boldsymbol{R}_c^e}{\mathrm{D}\bar{\boldsymbol{v}}} = \int_\Omega (\boldsymbol{N}^T\boldsymbol{N} \odot \boldsymbol{I})\frac{1}{\Delta t} + (\boldsymbol{N}^T \odot \boldsymbol{I})\nabla(\bar{\boldsymbol{v}} + \boldsymbol{v}')(\boldsymbol{N} \odot \boldsymbol{I}) + (\boldsymbol{N}^T(\bar{\boldsymbol{v}} + \boldsymbol{v}')^T \boldsymbol{J}^{-T}(\boldsymbol{DN})^T \odot \boldsymbol{I})$$
$$(81) \quad + 2\nu_f((\boldsymbol{DN})\boldsymbol{J}^{-1}\boldsymbol{J}^{-T}(\boldsymbol{DN})^T \odot \boldsymbol{I})\,\mathrm{d}\Omega$$

$$(82) \quad \frac{\mathrm{D}\boldsymbol{R}_p^e}{\mathrm{D}\bar{\boldsymbol{v}}} = \int_\Omega \boldsymbol{N}^T \mathrm{vec}[\boldsymbol{J}^{-T}(\boldsymbol{DN})^T]^T\,\mathrm{d}\Omega$$

$$(83) \quad \frac{\mathrm{D}\boldsymbol{R}_c^e}{\mathrm{D}\boldsymbol{v}'} = \int_\Omega (\boldsymbol{N}^T \odot \boldsymbol{I})b^e \nabla(\bar{\boldsymbol{v}} + \boldsymbol{v}') + (\boldsymbol{N}^T \nabla b^{eT}(\bar{\boldsymbol{v}} + \boldsymbol{v}') \odot \boldsymbol{I}) + 2\nu_f((\boldsymbol{DN})\boldsymbol{J}^{-1}\nabla b^e \odot \boldsymbol{I})\,\mathrm{d}\Omega$$

$$(84) \quad \frac{\mathrm{D}\boldsymbol{R}_p^e}{\mathrm{D}\boldsymbol{v}'} = \int_\Omega \boldsymbol{N}^T \nabla b^{eT}\,\mathrm{d}\Omega$$

$$\text{(85)} \quad \frac{\mathrm{D}\boldsymbol{R}_c^e}{\mathrm{D}p} = \int_\Omega -\mathrm{vec}[\boldsymbol{J}^{-T}(\boldsymbol{D}\boldsymbol{N})^T]\boldsymbol{N}\,\mathrm{d}\Omega$$

$$\begin{aligned}\text{(86)} \quad \frac{\mathrm{D}\boldsymbol{R}_f^e}{\mathrm{D}\bar{\boldsymbol{v}}} = \int_{\Omega_e} & b^e\boldsymbol{I}(\boldsymbol{N}\odot\boldsymbol{I})\frac{1}{\Delta t} + b^e\boldsymbol{I}\nabla(\bar{\boldsymbol{v}}+\boldsymbol{v}')(\boldsymbol{N}\odot\boldsymbol{I}) + b^e\boldsymbol{I}((\bar{\boldsymbol{v}}+\boldsymbol{v}')^T\boldsymbol{J}^{-T}(\boldsymbol{D}\boldsymbol{N})^T\odot\boldsymbol{I}) \\ & + 2\nu_f(\nabla b^{eT}\boldsymbol{J}^{-T}(\boldsymbol{D}\boldsymbol{N})^T\odot\boldsymbol{I})\,\mathrm{d}\Omega\end{aligned}$$

$$\text{(87)} \quad \frac{\mathrm{D}\boldsymbol{R}_f^e}{\mathrm{D}\boldsymbol{v}'} = \int_{\Omega_e} b^{e2}\nabla(\bar{\boldsymbol{v}}+\boldsymbol{v}') + (\nabla b^{eT}(\bar{\boldsymbol{v}}+\boldsymbol{v}'))b^e\boldsymbol{I} + 2\nu_f(\nabla b^{eT}\nabla b^e\odot\boldsymbol{I})\,\mathrm{d}\Omega$$

$$\text{(88)} \quad \frac{\mathrm{D}\boldsymbol{R}_f^e}{\mathrm{D}p} = \int_{\Omega_e} -(\nabla b^{eT}\odot\boldsymbol{I})\mathrm{vec}[\boldsymbol{I}]\boldsymbol{N}\,\mathrm{d}\Omega$$

## ACKNOWLEDGMENTS

The research reported herein was supported by the Computational Science and Engineering Fellowship (D. Z. Turner) and The Department of Energy (K. B. Nakshatrala) through a SciDAC-2 project (Grant No. DOE DE-FC02-07ER64323). This support is gratefully acknowledged. The opinions expressed in this paper are those of the authors and do not necessarily reflect that of the sponsor.

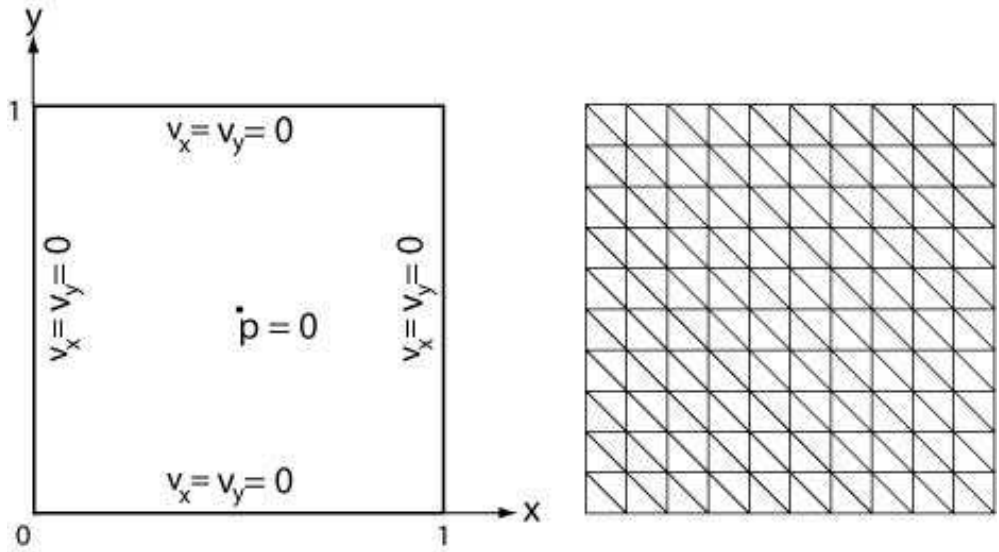

Figure 1. Body force driven cavity: problem statement and boundary conditions (left) computational mesh (right).

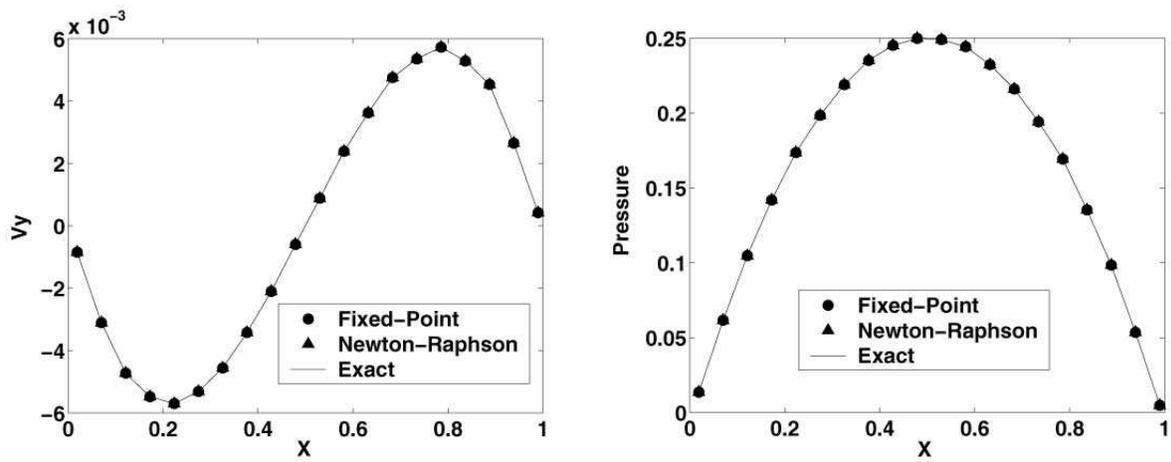

Figure 2. Body force driven cavity: pressure and velocity in the y-direction as measured along the x-axis at $y = 0.5$.

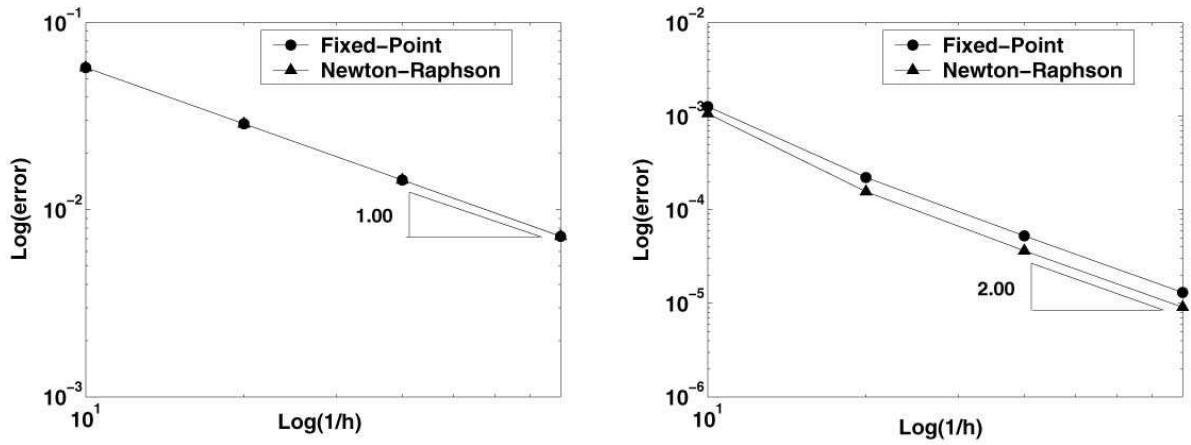

FIGURE 3. Body force driven cavity: $H^1$-semi norm of the pressure vs. element size for $Re = 400$ (left) $L^2$ norm of the velocity vs. element size (right).

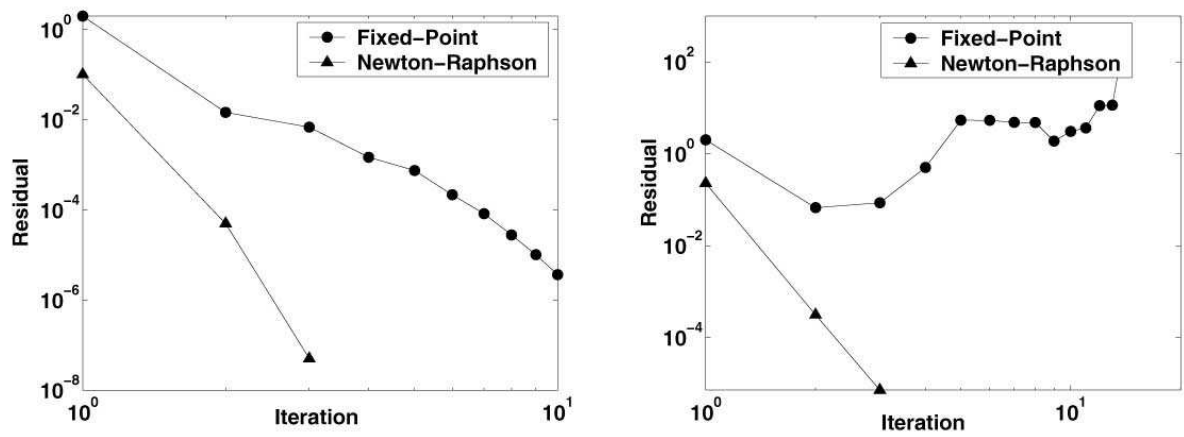

FIGURE 4. Body force driven cavity: convergence of the residual for $Re = 400$ (left) $Re = 5,000$ (right).

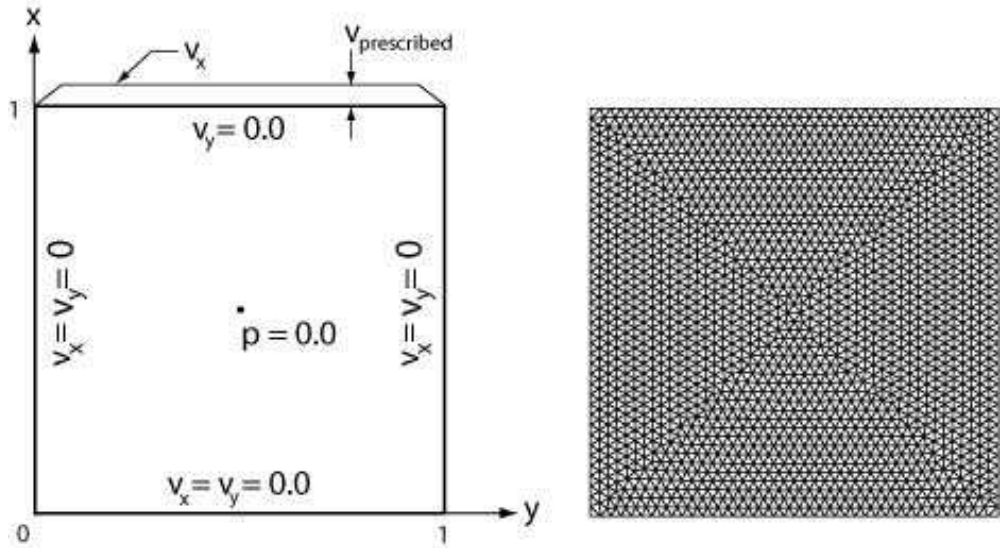

Figure 5. Lid driven cavity: problem statement and boundary conditions (left) computational mesh (right). The non-leaky cavity approach is used here which resolves the discontinuity at the upper two corners of the domain by assuming that the corners belong to the vertical walls.

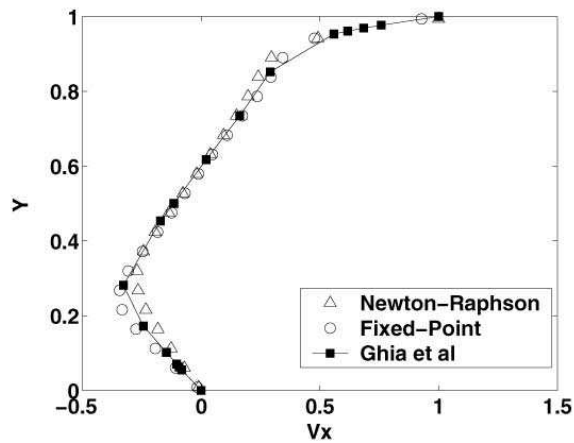

Figure 6. Lid driven cavity: velocity in the x-direction as measured along the y-axis at $x = 0.5$.

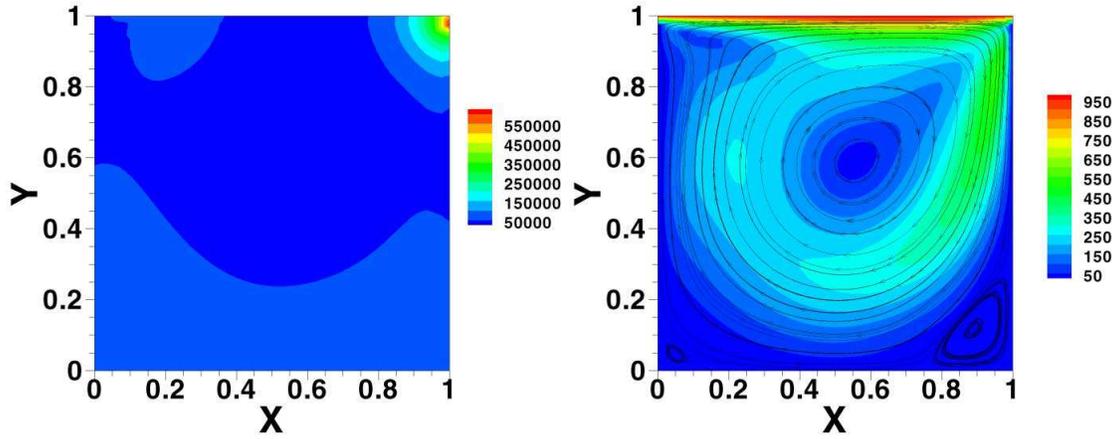

Figure 7. Lid driven cavity: velocity streamtraces over velocity magnitude contours for $Re = 400$ using the consistent Newton-Raphson approach.

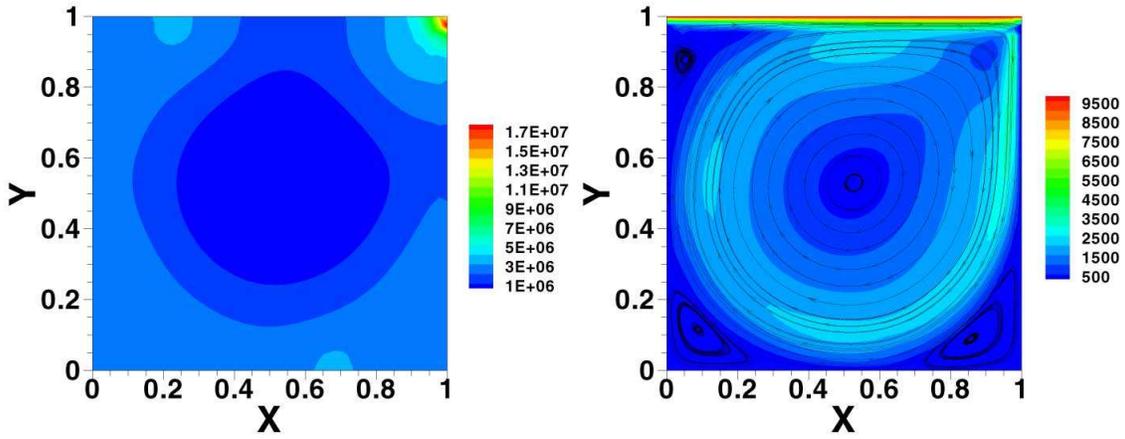

Figure 8. Lid driven cavity: velocity streamtraces over velocity magnitude contours for $Re = 5,000$ using the consistent Newton-Raphson approach.

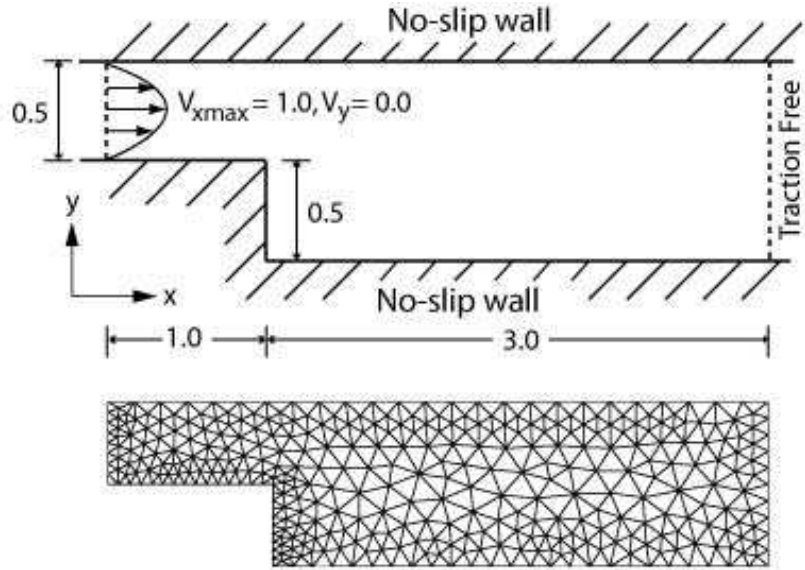

Figure 9. Backward facing step: problem statement and boundary conditions (top) computational mesh (bottom).

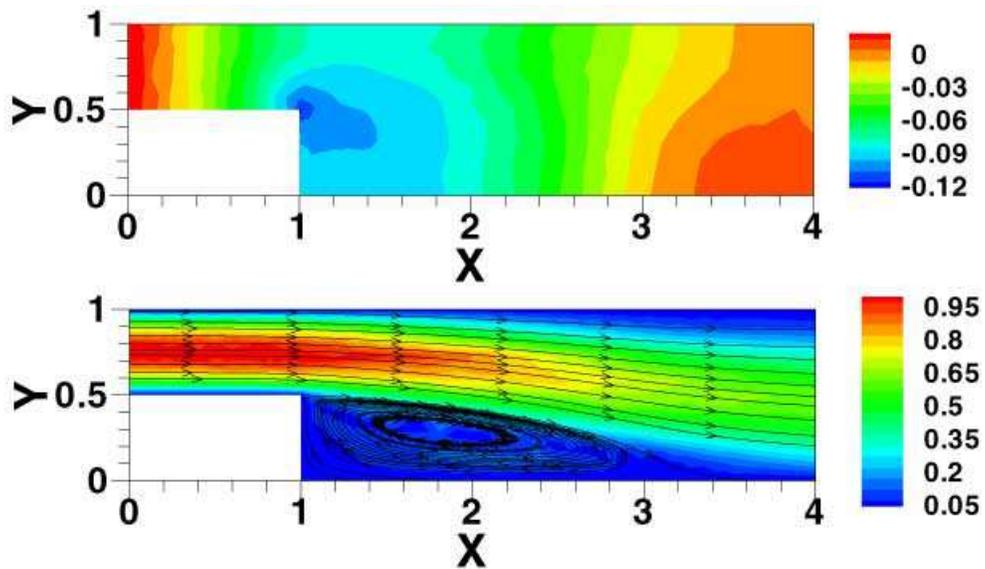

Figure 10. Backward facing step: using the consistent Newton-Raphson method for $Re = 150$, velocity streamtraces over velocity magnitude contours (top) pressure contours (bottom).

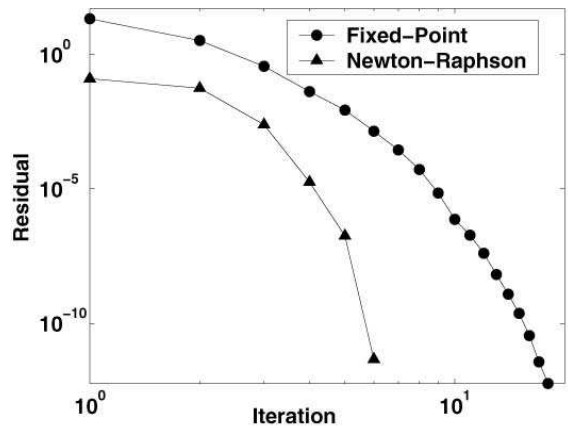

Figure 11. Backward facing step: convergence of the residual for $Re = 15$.

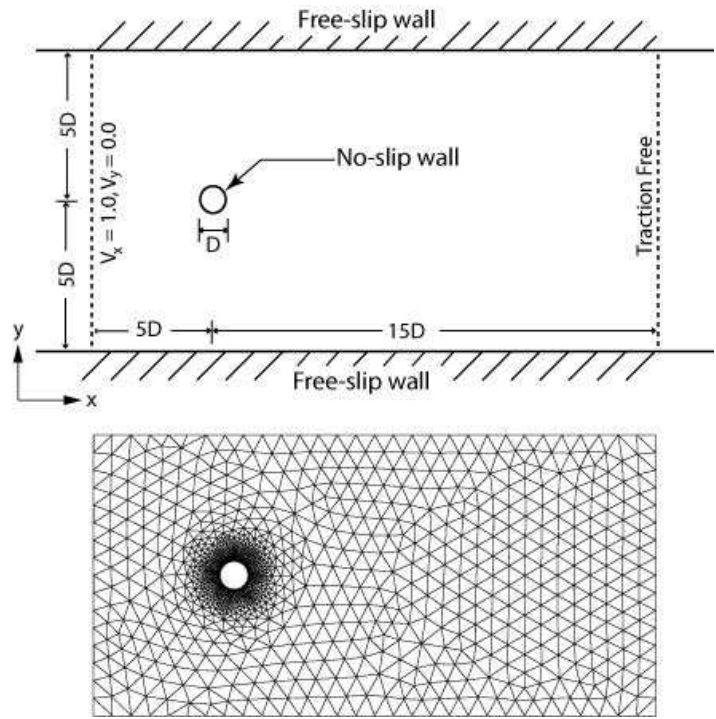

Figure 12. Flow past a circular cylinder: problem statement and boundary conditions (top) computational mesh (bottom).

Correspondence to: Daniel Z. Turner, Department of Civil and Environmental Engineering, 2103 Newmark Laboratory, University of Illinois at Urbana-Champaign, Urbana, Illinois - 61801.

*E-mail address*: dzturne1@illinois.edu

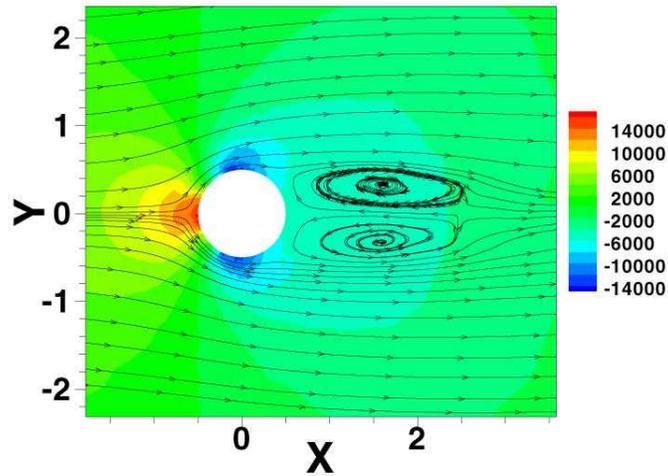

FIGURE 13. Flow past a circular cylinder: streamtraces shown over pressure contours for $Re = 75$.

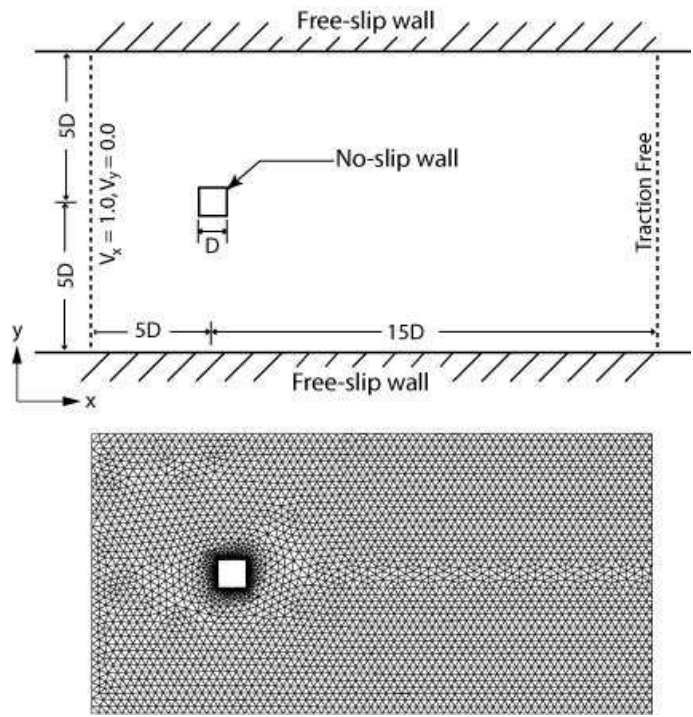

FIGURE 14. Flow past square cylinder: problem statement and boundary conditions (top) computational mesh (bottom).

Dr. Kalyana Babu Naskshatrala, Department of Civil and Environmental Engineering, 2524 Hydrosystems Laboratory, University of Illinois at Urbana-Champaign, Urbana, Illinois - 61801.

*E-mail address*: `nakshatr@illinois.edu`

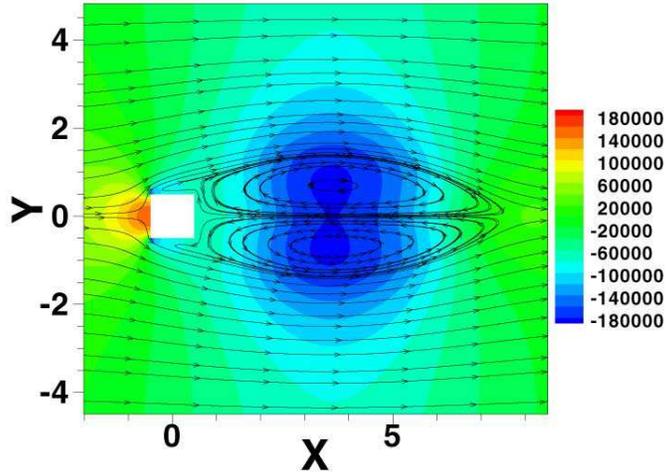

FIGURE 15. Flow past a square cylinder: streamtraces shown over pressure contours for $Re = 75$.

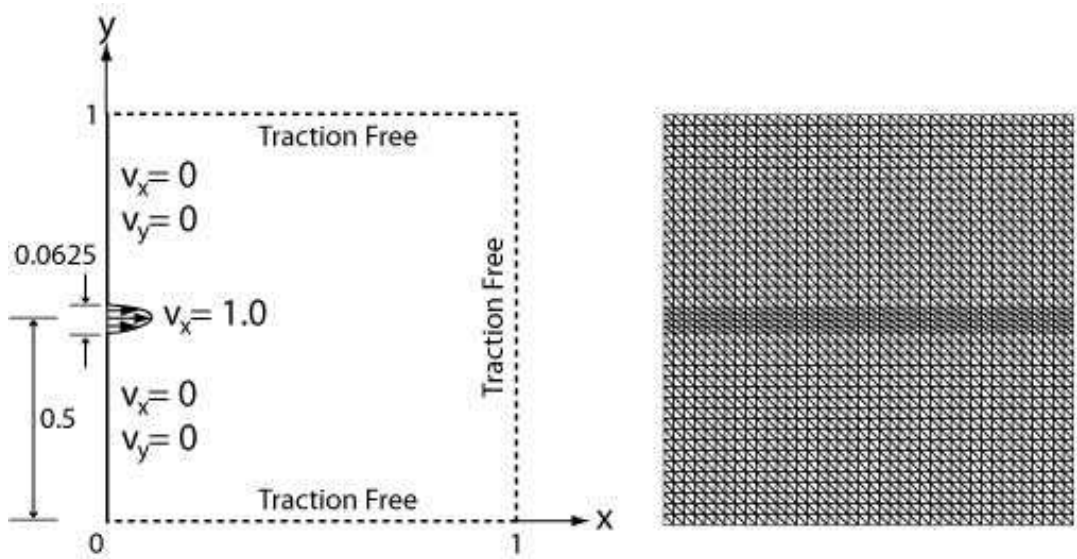

FIGURE 16. Jet flow: problem statement and boundary conditions (left) computational mesh (right).

Professor Keith D Hjelmstad, Department of Civil and Environmental Engineering, 3129e Newmark Laboratory, University of Illinois at Urbana-Champaign, Urbana, Illinois - 61801.

*E-mail address*: `kdh@illinois.edu`

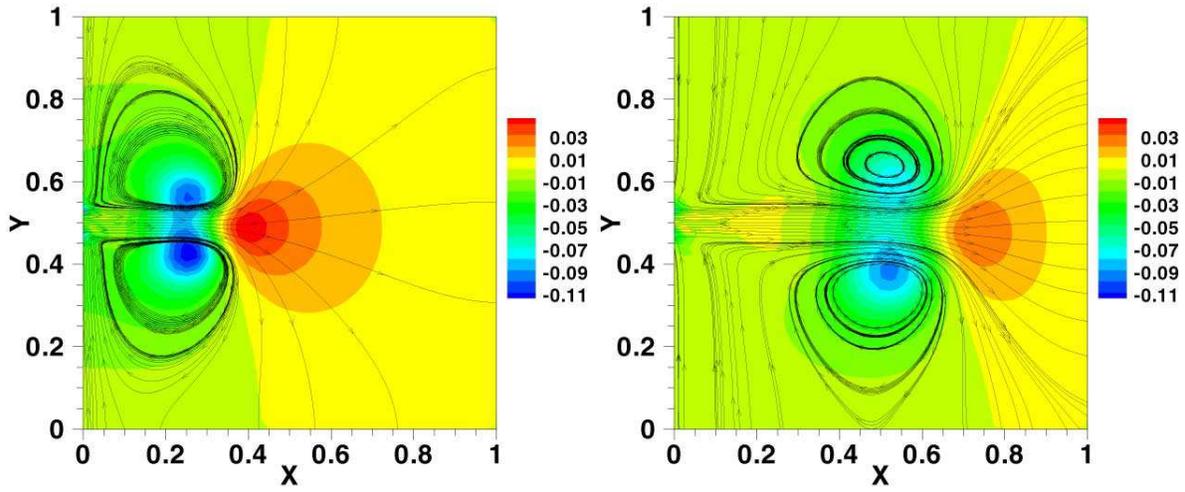

FIGURE 17. Jet flow: snapshot of velocity streamtraces over the pressure contours at $t = 1.0$s and $t = 2.24$s.